\documentclass[12pt]{amsart}
\usepackage{etex}
\usepackage{amscd,amssymb,hyperref}
\usepackage[PostScript=dvips]{diagrams}
\input xy
\usepackage{graphics}
\usepackage{epsfig}
\usepackage{picinpar}
\usepackage{amsmath,amsfonts,amsthm}
\usepackage{color}
\usepackage{pstricks,pst-node,pst-text,pst-3d}

\usepackage{verbatim}

\usepackage{pgfplots}

\pgfplotsset{every axis/.append style={line width=1pt}}

\usepackage{color}

\definecolor{BlueGreen}{RGB}{49,152,255}
\definecolor{DarkGreen}{RGB}{0,100,0}
\definecolor{DeepPink}{RGB}{255,20,147}
\definecolor{VioletRed}{RGB}{208,32,144}

\newcounter{EmptyBlock}[section]
{\par}%

\newtheorem{thm}{Theorem}[section]

\newtheorem{cor}[thm]{Corollary}

\def\square{\vbox{
      \hrule height 0.4pt
      \hbox{\vrule width 0.4pt height 5.5pt \kern 5.5pt \vrule width 0.4pt}
      \hrule height 0.4pt}}

\def\ch\mathrm{c h}

\newcommand{\N}{\mathbb{N}}

\numberwithin{equation}{section}

\title{On the generating function 
and growth 
of the positive  singular braid monoid 
}

\author{A.~L.~Anisimov%
}
\address{Nosov Magnitogorsk State Technical University, Magnitogorsk 455000, Lenin prospect, 38,
Russia}
\email{aanisimov@yandex.ru}
\author{G.~A.~Kameneva%
}
\address{Nosov Magnitogorsk State Technical University, Magnitogorsk 455000, Lenin prospect, 38, 
Russia}
\email{kameneva\_galina@mail.ru}

\author{V.~V.~Vershinin%
}

\address{Institut Montpelliérain Alexander Grothendieck, 
                                     Universit\'e de Montpellier,
Place Eug\`ene Bataillon,
34095 Montpellier cedex 5, France}
\email{vladimir.verchinine@umontpellier.fr}
\address{Sobolev Institute of Mathematics, Novosibirsk 630090,
Russia }

\numberwithin{figure}{section}



\begin{document}
\begin{abstract}
We give an elementary  proof of the fact that the generating function of the positive singular braid monoid 
is rational and we give the 
exact formula for such monoid on three strands.	
\end{abstract}

\maketitle

\section{Introduction}

Braid groups were defined by E.~Artin \cite{Artin} in 1925 and since then they were studied intensively.
They were also generalized in several directions. E. Brieskorn \cite {Bri1} defined what are called now 
Artin groups of finite type. 

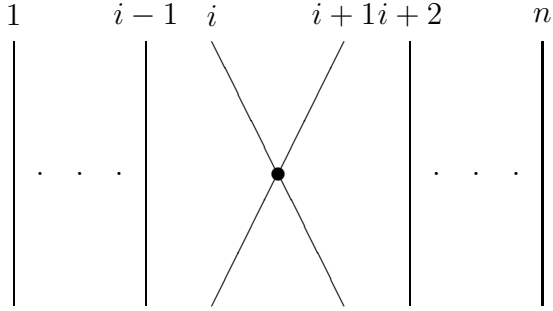
\begin{figure}
\begin{picture}(0,130)(0,-10) 
	\put(0,50){\circle*{5}} \put(-100,100){\line(0,-1){100}}
	\put(-50,100){\line(0,-1){100}} \put(-25,100){\line(1,-2){50}}
	\put(25,100){\line(-1,-2){50}} \put(50,100){\line(0,-1){100}}
	\put(100,100){\line(0,-1){100}}
	\put(-100,110){\makebox(0,0)[cc]{$1$}}
	\put(-50,110){\makebox(0,0)[cc]{$i-1$}}
	\put(-25,110){\makebox(0,0)[cc]{$i$}}
	\put(25,110){\makebox(0,0)[cc]{$i+1$}}
	\put(50,110){\makebox(0,0)[cc]{$i+2$}}
	\put(100,110){\makebox(0,0)[cc]{$n$}}
	\put(-75,50){\makebox(0,0)[cc]{.\quad.\quad.}}
	\put(75,50){\makebox(0,0)[cc]{.\quad.\quad.}}
\end{picture}
	\caption{The singular braid generator $x_i$}\label{fi:singen}
\end{figure}

The {\it Baez--Birman monoid} $SB_n$  or {\it singular braid monoid} 
\cite{Bae}, \cite{Bir2}
is defined as the
monoid with generators $\sigma_i,\sigma_i^{-1},x_i$, $i=1,\dots,n-1,$ and
relations 
\begin{eqnarray*}
	&\sigma_i\sigma_j=\sigma_j\sigma_i, \ \text {if} \ \ |i-j| >1,\\
	&x_ix_j=x_jx_i, \ \text {if} \ \ |i-j| >1,\\ &x_i\sigma_j=\sigma_jx_i, \
	\text {if} \ \ |i-j| \not=1,\\ &\sigma_i \sigma_{i+1} \sigma_i = \sigma_{i+1} \sigma_i
	\sigma_{i+1},\\ &\sigma_i \sigma_{i+1} x_i = x_{i+1} \sigma_i \sigma_{i+1},\\ &\sigma_{i+1} \sigma_i
	x_{i+1} = x_i \sigma_{i+1} \sigma_i,\\ &\sigma_i\sigma_i^{-1}=\sigma_i^{-1}\sigma_i =1.
\end{eqnarray*}
In pictures $\sigma_i$ corresponds to canonical generator of the braid
group and $x_i$ represents an intersection
of the $i$th and $(i+1)$st strand as in the 
Figure~\ref{fi:singen}. 

Let $G$ be a finitely presented group or monoid, and let $a_k$ be the number of its elements of the length k.
One of the aspects of study of groups and monoids  is to understand the behaviour of  the 
sequence $(a_k)$, $k\in \N$.  

Let $B(W)_n^+$ is the positive
Artin monoid of finite type (no inverses of generators), $W$ is the corresponding finite Coxeter group.
Let $a_k^{[n]}$ be the number of elements of the length $k$ in such monoid 
and 
\begin{align*}
	\phi^{[n]}(t)=\sum_{k\geq 0} a^{[n]}_k t^k
\end{align*}
be the corresponding generating function.

P.~Deligne proved that $ \phi^{[n]}(t)$ is a rational function. Independently Xu
proved this for the classical positive braid monoids and gave the following expression of this
function for $B_3^+$
\begin{align*}
	\frac 1 {(1-x) (1-x-x^2)}. 
\end{align*}
So $a_k=F_{k+3}-1$, where $F_i$ are the Fibonacci numbers. 
Hence the growth rate of the classical braid monoid on three strands, is equal
to $\frac{1+ \sqrt 5} 2 \approx 1.618$.
Concrete formulas for the generating functions for Artin monoids were obtained by 
K.~Saito etal \cite {S, FFST}. B.~Berceanu and Z.~Iqbal proved in  \cite{BI}  that 4 is a
universal upper bound for the growth of the positive Artin monoids of finite type.

Let $SB_n^+$  the positive singular braid 
monoid, i.e. given by the  generators $\sigma_i, x_i$, $i=1,\dots,n-1,$ and the positive part 
of singular braid 
relations. 
Let $b_k^{[n]}$ be its number of elements of the length $k$ in such monoid 
and 
\begin{align*}
	f^{[n]}(t)=\sum_{k\geq 0} b^{[n]}_k t^k
\end{align*}
be the corresponding generating function. The monoids $SB_n^+$ have been deeply studied,
including growth functions by R. Corran \cite{Cor, Cor2}. In thie present work we propose an elementary approach.

\section{Statements of the main results}

Here we prove the following theorems.
\begin{thm}\label{theorem1}
The  generating function for the positive singular braid monoid $SB_n^+$ is rational.
\end{thm}
The monoid $SB_3^+$ is given by four generators $\sigma_1, \sigma_2, x_1, x_2$ and following relations.
\begin{align*}
	\left\{\begin{array}{l} \sigma_2\sigma_1\sigma_2 = \sigma_1\sigma_2\sigma_1, \\ 
		x_2\sigma_1\sigma_2 = \sigma_1\sigma_2 x_1, \\
		x_1\sigma_2\sigma_1 = \sigma_2\sigma_1 x_2, \\  
		x_1 \sigma_1 = \sigma_1 x_1, \\
		x_2 \sigma_2 = \sigma_2 x_2.
	\end{array}\right.   
\end{align*}
\begin{thm}\label{theorem2}
The generating function for the monoid $SB_3^{+}$ is given by the formula
\begin{equation} \label{sysx}
	f^{[3]} =	\frac{1} {(1-t)(1-3t-t^2+ 2t^3)} .
\end{equation}
\end{thm}

\begin{cor}\label{cor}
	The sequence of numbers $b_n$, $n\geq 0$, satisfies the condition
	$$b_0=1
	$$ 
	and the following recurrence
relations 	\begin{equation*}
		b_k = 4 b_{k-1} - 2 b_{k-2} - 3 b_{k-3}  + 2 b_{k-4}  \ \  n \geqslant 1,
	\end{equation*}
or
\begin{equation*}
b_k = 3 b_{k-1} + b_{k-2} - 2 b_{k-3} + 1, \;  \ \ n \geqslant 1,
\end{equation*}
where all $b_k$ with $k\ <0$ are equal to $0$.
	\end{cor}

In particular we have
$$
b_1= 4, \ b_2= 14, \ b_3=45, \  b_4= 142, \ b_5 = 444, \dots .
$$
For $k=2$ the are the following elements $ \sigma_1^2, \sigma_2 \sigma_1, \sigma_1\sigma_2, \sigma_2^2, 
x_1 \sigma_1,  x_1 \sigma_2, $ $x_2\sigma_2,  x_2 \sigma_2, \sigma_2 x_1, \sigma_1 x_2, x_1^2, x_2^2, x_2x_1, x_1x_2$.

\section{Proofs}

We use the right greedy form of Adyan \cite{Ad}, El-Rifai \& Morton \cite{EM} and Thurston etal \cite{E_Th}
which was adapted for singular braids
in \cite{Cor} and \cite{Ver}. We also use the constructions and proofs of \cite{Xu}. 

We denote by $\{g_i  \},$  the  ordered set of all divisers of the fundamental word $\Delta_n$. 
Its cardinality is $n!$. We denote by $g_{n!+k}$   the generators  $x_k$, $k=1, \dots , n-1$.
The length of $g_i$ is denoted by $l_i$.
We denote by $M_n$ the positive singular braid monoid $SB_n^+$, and let $m=n!+n-1$. 
Let $f^{[n]}$ be the generating function 
of the  monoid $M_n$ with respect to the canonical generators.   Let  $M_n^j$ be the subset of 
elements of $M_n$ whose maximal (among all $g_j$) right diviser  is $g_j$, j$=1, \dots , m$.
 Let $f_j$ be the generating function 
 of $M_n^j$, then
 $$
  f^{[n]}(t) =\sum_{j=1}^m f_j(t).
  $$ 
$M_n^1=\{1\}$  and so $f_1(t)=1$. 
We define the set of predecessors ${\rm Pred} (g_i)$ of $g_i$  as the set of those $g_j$ that can 
precede $g_i$
in the right greedy normal form. Let ${\rm Pred}_\Delta (g_i)$, $i\leq n!$ be the set of predecessors of 
$g_i$ for the classical positive braid monoid.  For $g_i$, $i\leq n!$ we have
$$ {\rm Pred}_(g_i) = {\rm Pred}_\Delta (g_i) \cup \{x_!, \dots, x_{n-1} \}.
$$
For $x_k$ the set of predecessors consists of all $x_j$ and those $g_i)$, $i\leq n!$ that do not 
have $\sigma_k$ and $\sigma_k\sigma_{k-1}$ as right divisors.
Table 1 describes the predecessors in the case of the monoid on 3 strands.
We define the numbers $\epsilon_{i,j}$ as follows
$$
\epsilon_{i,j} = \begin{cases} 
		0,  & \text{if }  \ g_j\not\in {\rm Pred} (g_i),\\
                                  1,  & \text{if } g_j\in {\rm Pred} (g_i).
\end{cases}
$$
\subsection{Proof of Theorem~1}

From the definitions of $f_i$ we have the following system of recursive linear equations
\begin{equation} \label{sys}
\frac{f_i}{t^{l_i}} = 1+\sum_{j=2}^m \epsilon_{i,j} f_j , \ \ i= 2, \dots, m.
\end{equation}
The element $\Delta$ can not precede any other $g_i$, but it can follow any other $g_i$.
So $\epsilon_{i,n!}=0$ for $i\not= n!$, $\epsilon_{n!,j}=1$ for all $j$.  The length of $\Delta$ we denote by $l$.
Hence we can rewrite the system \ref{sys} in the form
\begin{equation} \label{sys2}
\begin{cases} 	\frac{f_i}{t^{l_i}} = 1+&\sum_{j=2}^{n!-1} \epsilon_{i,j} f_j  + \sum_{j=n!+1}^{m} \epsilon_{i,j} f_j , \
	i= 2, \dots, n!-1, n!+1, \dots, m,\\
		\frac{f_{n!}}{t^l} = 1+&\sum_{j=2}^{m} f_j. 
		\end{cases}
\end{equation}
The first $m-2$ equations form a system with $m-2$ variables. Its determinant is a non-zero rational function on $t$. 
So, each $f_i$, $i\not=n!$ is a rational function.
Then we express $f_{n!}$ according to (\ref{sys2}.)
\begin{equation*}
	f_{n!}=  \frac {t^l(1+ \sum_{j=2}^{n!-1} f_j  + \sum_{j=n!+1}^{m}  f_j )}{1-t^l}
\end{equation*}
We get an expression 
\begin{equation} \label{sysy} 
	f^{[n]}=  \frac {1+ \sum_{j=2}^{n!-1} f_j  + \sum_{j=n!+1}^{m}  f_j }{1-t^l}
\end{equation}
$\square$

\subsection{Proof of Theorem~2}

Let $n=3$. 
According to Table 1 
the system for $f_2, f_3, f_4, f_5, f_7, f_8$
is the following.
\begin{equation*} \label{sys3}
	\begin{cases} 	(1-t^{-1})f_2 + f_5 +f_7+f_8 =&-1\\
		(1-t^{-1})f_3 + f_4 +f_7+f_8 =&-1\\
		f_2 -t^{-2} f_4 +f_5+ f_7+f_8 =&-1\\
		f_3 +f_4 -t^{-2} f_5 +f_7+f_8 =&-1\\
			f_3 +(1-t^{-1}) f_7 +f_8 =&-1\\
			f_2 +f_7 +(1-t^{-1})f_8 =&-1\\
			\end{cases}
\end{equation*}

\hypertarget{handtab:1}{}
\phantomsection
\addcontentsline{toc}{subsection}{Table 1. Information data for $SB_3^{+}$.}
\subsection*{Table 1. Information data for $SB_3^{+}$.}
\begin{center}
	\begin{tabular}{|c||c|c|c|}                                       \hline
		$i$   & Braid rep. of $g_i$  & $l_i$  & $Pred(g_i)$            \\ \hline
		$2$   & $\sigma_1$           & $1$    & $g_2, g_5, g_7, g_8$   \\ \hline
		$3$   & $\sigma_2$           & $1$    & $g_3, g_4, g_7, g_8$   \\ \hline
		$4$   & $\sigma_1\sigma_2$   & $2$    & $g_2, g_5, g_7, g_8$   \\ \hline
		$5$   & $\sigma_2\sigma_1$   & $2$    & $g_3, g_4, g_7, g_8$   \\ \hline
		$7$   & $x_1$                & $1$    & $g_3, g_7, g_8$        \\ \hline
		$8$   & $x_2$                & $1$    & $g_2, g_7, g_8$        \\ \hline
	\end{tabular}
\end{center}

There is the automorphism of $SB_3^+$ given by
\begin{equation*} 
	\begin{cases} 	\sigma_1\to &\sigma_2\\
		\sigma_2\to &\sigma_1\\
			x_1\to &x_2\\
			x_2\to &x_1\\
			\end{cases}
\end{equation*}
So, $f_2=f_3$, $f_4=f_5$, $f_7=f_8$ and we have the following system
\begin{equation*} \label{sys5}
	\begin{cases} 	(1-t^{-1})f_2 + f_4 +2f_7 =&-1\\
				f_2 +(1-t^{-2}) f_4  +2f_7 =&-1\\
				f_2 +(2-t^{-1}) f_7  =&-1\\
		\end{cases}
\end{equation*}
From the first two equations we get $f_4=tf_2$ and from the third one that
\begin{equation*} \label{sys6}
		 f_7  =\frac t {1-2t} f_2+\frac t {1-2t}\\
\end{equation*}
Finally we get the equation for $f_2$:
\begin{equation*} \label{sys7}
	 	(1-t^{-1})f_2 + tf_2+2(\frac t {1-2t} f_2+\frac t {1-2t})=-1
		\end{equation*}
and the expressions for $f_2$ and $f_7$
	\begin{equation*} \label{sys3}
		f_2 =\frac t {1-3t-t^2+ 2t^3} \\
		\end{equation*}
	\begin{equation*} \label{sys3}
	f_7 =\frac {t^2} {(1-2t)(1-3t-t^2+ 2t^3)}  +\frac t {1-2t}
\end{equation*}
We rewrite (\ref{sysy}) in this case
	\begin{equation*} 
	f^{[3]} =\frac{1 +f_2 +f_3 +f_4+f_5 +f_7+f_8} {(1-t^3)} 
\end{equation*}
Using formulas for $f_i$  we get an expression for $f^{[3]}$
	\begin{equation*} 
		f^{[3]} =\frac{1-t-t^2-2t^3} {(1-t^3)(1-2t)(1-3t-t^2+ 2t^3)} =\\
		\frac{1} {(1-t)(1-3t-t^2+ 2t^3)} .
\end{equation*}
$\square$

\section{Growth}

Let us consider the equation
\begin{equation*}
2t^3-t^2-3t+1=0.
\end{equation*}
Its depressed form is as follows
\begin{equation*}
	y^3-\frac {19} {12} y+\frac{13}{54} =  0.
\end{equation*}
The value of 
\begin{equation*}
	\Delta= \frac{q^2} 4 + \frac {p^3} {27} = \frac 1 4 \left({\frac{13}{54}}\right)^2  - \frac 1 {27} \left(\frac{19}{12}\right)^3
\end{equation*}
is negative. So, the equation has three real roots. 
Let us denote these roots by $r_1$, $r_2$ and $r_3$ in ascending order. Their approximative values are as follows
\begin{align*}
	r_1 &\approx -1.161702138, \\
	r_2 &\approx 0.3210368161, \\
	r_3 &\approx 1.340665322.
\end{align*}

They can be also expressed by the trigonometric formulas. For example for $r_2$ : 
\begin{equation}\label{r}
	r_2 = \frac{\sqrt{19}\cos{\left(\frac{\arccos{\left(\frac{26\sqrt{19}}{361}\right)}}{3}+\frac{\pi}{3}\right)}}{3} + \frac16 \approx 0.3210368161, \\
\end{equation}
Let us rewrite the function (\ref{sysx}) in the form
\begin{equation} \label{sysx2}
		f^{[3]} =
		\frac{1} {(1-t)(1-3t-t^2+ 2t^3)}= \\
\frac{1} {2(1-t)(r_1-t)(r_2-t)(r_3-t)}	.
\end{equation}
and then we express (\ref{sysx2}) as a sum with 
coefficients $a_i$, $i=0,1,2,3$.
\begin{equation*} \label{sysx3}
		f^{[3]} =
\frac{a_0} {(1-t)}+ \frac{a_1} {(r_1-t)}+	\frac{a_2} {(r_2-t)}
+ \frac{a_3} {(r_3-t)}	.
\end{equation*}
After calculations we obtain:
\begin{equation*} \label{sysx31}
	a_0 = -1, 
\end{equation*}
\begin{equation*}
	a_1 = \frac12 \frac{2r_1^2+r_1-2}{(r_1-r_2)(r_1-r_3)} \approx -0.06233879045, 
\end{equation*}
\begin{equation*}
	a_2 = \frac12 \frac{2r_2^2+r_2-2}{(r_2-r_1)(r_2-r_3)} \approx 0.4870988600, 
\end{equation*}
\begin{equation*}
	a_3 = -\frac12 \frac{2r_3^2+r_3-2}{(r_3-r_1)(r_3-r_2)} \approx 0.5752399310.
\end{equation*}
We have 
\begin{equation*} \label{sysx4}
\frac{a_j} {(r_j-t)}= a_j \sum_{i=0}^{\infty} r_j^{-(i+1)} t^i,	\
\ j= 0, 1,  2,  3.
\end{equation*}
As $ | \, r_j| > 1$ in the study of the growth rate  the only essential series is 
\begin{equation*} \label{sysx5}
\frac{a_2} {(r_j-t)}= a_2 \sum_{i=0}^{\infty} r_2^{-(i+1)} t^i	.
\end{equation*}
So the growth rate of the $SB_3^+$ is equal to $r_2^{-1}$ where $r_2$ is given by the formula  (\ref{r}). Approximatively it is equal to 3.11, almost two times greater than that of the growth rate of the classical braid monoid on three strands, which is equal to $\frac{1+ \sqrt 5} 2 \approx 1.618$. This looks quite
reasonable because of two new generators.

\end{document}